\documentclass[12pt,twoside]{article}
\usepackage{amsmath,amsfonts,amssymb,theorem}
\textheight215truemm \textwidth 138truemm  \date{}
 \newcommand\C{\mathbb C}
 \newcommand\N{\mathbb N}
 \newcommand\Q{\mathbb Q}
 \newcommand\R{\mathbb R}
 \newcommand\Z{\mathbb Z}
 \newcommand\A{\mathbb A}

 \newcommand{\iso}{\cong}
 \newcommand{\dd}{\mathrm d} 
 \newcommand\bdt{\mathbf{dt}} 
 \newcommand\bdx{\mathbf{dx}} 
 \newcommand\xbar{\overline x}

 \newcommand\al{\alpha}
 \newcommand\be{\beta}
 \newcommand\ga{\gamma}
 \newcommand\om{\omega}
 \newcommand\fie{\varphi}
 \newcommand\Ga{\Gamma}
 \newcommand\Om{\Omega}

 \newcommand{\cF}{\mathcal F}
 \newcommand{\Oh}{\mathcal O}
 \newcommand{\cR}{\mathcal R}
 \newcommand{\cS}{\mathcal S}
 \newcommand{\cV}{\mathcal V}
 \newcommand{\cW}{\mathcal W}
 \newcommand{\cY}{\mathcal Y}
 \newcommand{\cZ}{\mathcal Z}

 \newcommand{\fq}{\mathfrak q}
 \newcommand{\fU}{\mathcal U}
 \newcommand{\cX}{\mathcal X}

 \newcommand{\codim}{\operatorname{codim}}
 \newcommand{\Ord}{\operatorname{Ord}}
 \newcommand{\SL}{\operatorname{SL}}
 \newcommand{\Spec}{\operatorname{Spec}}

 \newtheorem{prop}{Proposition}[section]
 \newtheorem{thm}[prop]{Theorem}

 \theorembodyfont{\rmfamily}
 \newtheorem{rem}[prop]{Remark}
 \newtheorem{dfn}[prop]{Definition}

 \newenvironment{pf}{\paragraph{Proof}}{\par\medskip}
 \newenvironment{pfof}[1]{\paragraph{Proof of #1}}{\par\medskip}
 \newcommand{\qed}{\ifhmode\unskip\nobreak\fi\quad $\square$}

 \title{Birational Calabi--Yau $n$-folds\\ have equal Betti numbers}
 \author{Victor V. Batyrev}

\begin{document}

\maketitle

 \begin{abstract}
Let $X$ and $Y$ be two smooth projective $n$-dimensional algebraic
varieties $X$ and $Y$ over $\C$ with trivial canonical line bundles. We
use methods of $p$-adic analysis on algebraic varieties over local
number fields to prove that if $X$ and $Y$ are birational, they have the
same Betti numbers.
 \end{abstract}

\section{Introduction}

The purpose of this note is to show how to use the elementary
theory of $p$-adic integrals on algebraic varieties to prove
cohomological properties of birational algebraic varieties over
$\C$. We prove the following theorem, which was used by
Beauville in his recent explanation of a Yau--Zaslow formula for
the number of rational curves on a K3 surface \cite{Beauville}
(see also \cite{FGS,YZ}):

 \begin{thm}
 \label{main-th}
Let $X$ and $Y$ be smooth $n$-dimensional irreducible
projective algebraic varieties over $\C$. Assume that the
canonical line bundles $\Om^n_X$ and $\Om^n_Y$ are trivial and
that $X$ and $Y$ are birational. Then $X$ and $Y$ have the same
Betti numbers, that is,
 \[
 H^i(X,\C) \iso H^i(Y,\C) \quad \text{for all $i\ge0$.}
 \]
 \end{thm}

Note that Theorem~\ref{main-th} is obvious for $n =1$, and for $n=2$,
it follows from the uniqueness of minimal models of surfaces with
$\kappa\ge0$, that is, from the property that any birational map between
two such minimal models extends to an isomorphism \cite{KMM}. Although
$n$-folds with $\kappa\ge0$ no longer have a unique minimal models for
$n\ge3$, Theorem~\ref{main-th} can be proved for $n=3$ using a result
of Kawamata (\cite{kawamata}, \S6): he showed that any two birational
minimal models of $3$-folds can be connected by a sequence of flops
(see also \cite{kollar}), and simple topological arguments show that if
two projective $3$-folds with at worst $\Q$-factorial terminal
singularities are birational via a flop, then their singular Betti
numbers are equal. Since one still knows very little about flops in
dimension $n\ge4$, it seems unlikely that a consideration of flops
could help to prove Theorem~\ref{main-th} in dimension $n\ge4$.
Moreover, Theorem~\ref{main-th} is false in general for projective
algebraic varieties with at worst $\Q$-factorial Gorenstein terminal
singularities of dimension $n\ge4$. For this reason, the condition in
Theorem~\ref{main-th} that $X$ and $Y$ are smooth becomes very
important in the case $n\ge4$. We remark that in the case of
holomorphic symplectic manifolds some stronger result is obtained in
\cite{H}.

\bigskip

\noindent
{\bf Acknowledgements:}
The author would like to thank Professors A. Beauville, B. Fantechi, 
L. G\"ottsche, K. Hulek, Y. Kawamata, M. Kontsevich, S. Mori, M. Reid 
and  D. van Straten for  fruitful discussions and stimulating e-mails. 

\section{Gauge forms and $p$-adic measures}

Let $F$ be a local number field, that is, a finite extension of the
$p$-adic field $\Q_p$ for some prime $p\in\Z$. Let $R\subset F$ be the
maximal compact subring, $\fq\subset R$ the maximal ideal,
$F_{\fq}=R/\fq$ the residue field with $|F_{\fq}|=q=p^r$.
We write
 \[
 N_{F/{\Q}_p}\colon F\to\Q_p 
 \]
for the standard norm, and $\| \cdot \|\colon F \to \R_{\ge 0}$
for the multiplicative $p$-adic norm
 \[
 a \mapsto \|a\|=p^{-\Ord(N_{F/\Q_p}(a))}.
 \]
Here $\Ord$ is the $p$-adic valuation.

 \begin{dfn} Let $\cX$ be an arbitrary flat reduced algebraic
$S$-scheme, where $S=\Spec R$. We denote by $\cX(R)$ the set of
$S$-morphisms $S\to\cX$ (or sections of $\cX\to S$). We call $\cX(R)$
the set of $R$-{\em integral points} in $\cX$. The set of sections of
the morphism $\cX\times_S\Spec F\to\Spec F$ is denoted by $\cX(F)$ and
called the set of $F$-{\em rational points} in $\cX$.
 \end{dfn}

 \begin{rem}\label{point}
 \begin{enumerate}
 \renewcommand{\labelenumi}{(\roman{enumi})}
 \item If $\cX$ is an affine $S$-scheme, then one can
identify $\cX(R)$ with the subset
 \[
 \bigl\{x\in\cX(F) \bigm| 
 \text{$f(x)\in R$ for all $f\in\Ga(\cX,\Oh_{\cX})$}\bigr\}
 \subset \cX(F).
 \]

 \item If $\cX$ is a projective (or proper) $S$-scheme, then
$\cX(R)=\cX(F)$.
 \end{enumerate}
 \end{rem}

Now let $X$ be a smooth $n$-dimensional algebraic variety over
$F$. We assume that $X$ admits an extension $\cX$ to a regular
$S$-scheme. Denote by $\Om^n_X$ the canonical line bundle of
$X$ and by $\Om^n_{\cX/S}$ the relative dualizing sheaf on
$\cX$.

Recall the following definition introduced by Weil \cite{weil}:

 \begin{dfn} A global section $\om\in\Ga(\cX,\Om^n_{\cX/S})$ is
called a {\em gauge form} if it has no zeros in $\cX$. By
definition, a gauge form $\om$ defines an isomorphism
$\Oh_\cX\iso\Om^n_{\cX/S}$, sending $1$ to $\om$. Clearly, a
gauge form exists if and only if the line bundle $\Om^n_{\cX/S}$
is trivial.
 \end{dfn}

Weil observed that a gauge form $\om$ determines a canonical
$p$-adic measure $\dd\mu_\om$ on the locally compact $p$-adic
topological space $\cX(F)$ of $F$-rational points in $\cX$.
The $p$-adic measure $\dd\mu_\om$ is defined as follows:

Let $x\in\cX(F)$ be an $F$-point, $t_1, \dots, t_n$ local $p$-adic
analytic parameters at $x$. Then $t_1, \dots , t_n$ define a $p$-adic
homeomorphism $\theta \colon U \to\A^n(F)$ of an open subset
$\fU\subset\cX(F)$ containing $x$ with an open subset
$\theta(\fU)\subset\A^n(F)$. We stress that the subsets
$\fU\subset\cX(F)$ and $\theta(\fU)\subset\A^n(F)$ are considered to be
open in the $p$-adic topology, not in the Zariski topology. We write
 \[
\om=\theta^*\left(g \dd t_1 \wedge \cdots
\wedge \dd t_n\right),
 \]
where $g=g(t)$ is a $p$-adic analytic function on $\theta(\fU)$
having no zeros. Then the $p$-adic measure $\dd\mu_\om$ on $\fU$
is defined to be the pullback with respect to $\theta$ of the
$p$-adic measure $\|g(t)\|\bdt$ on $\theta(\fU)$, where $\bdt$ is the
standard $p$-adic Haar measure on $\A^n(F)$ normalized by the condition
 \[\int_{\A^n(R)} \bdt =1.
 \]
It is a standard exercise using the Jacobian to check that two
$p$-adic measures $\dd\mu_\om',\dd\mu_\om''$ constructed as above on
any two open subsets $\fU',\fU''\subset\cX(F)$ coincide on the
intersection $\fU'\cap\,\fU''$.

 \begin{dfn} The measure $\dd\mu_\om$ on $\cX(F)$ constructed above is
called the {\em Weil $p$-adic measure} associated with the gauge form
$\om$.
 \end{dfn}

 \begin{thm}[\cite{weil}, Theorem~2.2.5] Let $\cX$ be a regular
$S$-scheme, $\om$ a gauge form on $\cX$, and $\dd\mu_\om$ the
corresponding Weil $p$-adic measure on $\cX(F)$. Then
 \[
 \int_{\cX(R)} \dd\mu_\om=\frac{|\cX(F_\fq)|}{q^n},
 \]
where $\cX({F_\fq})$ is the set of closed points of $\cX$ over
the finite residue field $F_\fq$.
 \label{weil1}
 \end{thm}

 \begin{pf} Let
 \[
 \fie\colon \cX(R)\to\cX({F_\fq})
 \quad\text{given by}\quad
 x \mapsto \xbar\in\cX({F_\fq}) 
 \]
be the natural surjective mapping. The proof is based on the idea that
if $\xbar\in\cX({F_\fq})$ is a closed ${F_\fq}$-point of $\cX$ and
$g_1,\dots,g_n$ are generators of the maximal ideal of $\xbar$ in
$\Oh_{\cX,\xbar}$ modulo the ideal $\fq$, then the elements
$g_1,\dots,g_n$ define a $p$-adic analytic homeomorphism
 \[
 \ga\colon\fie^{-1}(\xbar)\to\A^n(\fq),
 \]
where $\fie^{-1}(\xbar)$ is the fiber of $\fie$ over $\xbar$ and
$\A^n(\fq)$ is the set of all $R$-integral points of $\A^n$ whose
coordinates belong to the ideal $\fq \subset R$. Moreover, the
$p$-adic norm of the Jacobian of $\ga$ is identically equal to
$1$ on the whole fiber $\fie^{-1}(\xbar)$. In order to see the latter
we remark that the elements define an \'etale morphism $g\colon V \to
\A^n$ of some Zariski open neighbourhood $V$ of $\xbar\in\cX$. Since
$\fie^{-1}(\xbar) \subset V(R)$ and $g^*(\dd t_1,\wedge \cdots \wedge
\dd t_n)=h \om$, where $h$ is invertible in $V$, we obtain that $h$ has
$p$-adic norm $1$ on $\fie^{-1}(\xbar)$. So, using the $p$-adic analytic
homeomorphism $\ga$, we obtain
 \[
 \int_{\fie^{-1}(\xbar)} \dd\mu_\om=\int_{\A^n(\fq)}\bdt
 =\frac{1}{q^n}
 \]
for each $\xbar\in\cX({F_\fq})$.
 \qed\end{pf}

Now we consider a slightly more general situation. We assume only that
$\cX$ is a regular scheme over ${S}$, but do not assume the existence
of a gauge form on $\cX$ (that is, of an isomorphism
$\Oh_\cX\iso\Om^n_{\cX/S}$). Nevertheless under these weaker
assumptions we can define a unique natural $p$-adic measure $\dd\mu$ at
least on the compact $\cX(R)\subset\cX(F)$ -- although possibly not on
the whole $p$-adic topological space $\cX(F)$!

Let $\fU_1,\dots,\fU_k$ be a finite covering of $\cX$ by Zariski open
$S$-subschemes such that the restriction of $\Om^n_{\cX/S}$ on each
$\fU_i$ is isomorphic to $\Oh_{\fU_i}$. Then each $\fU_i$ admits a gauge
form $\om_i$ and we define a $p$-adic measure $\dd\mu_i$ on each compact
$\fU_i(R)$ as the restriction of the Weil $p$-adic measure
$\dd\mu_{\om_i}$ associated with $\om_i$ on $\fU_i(F)$. We note that
the gauge forms $\om_i$ are defined uniquely up to elements
$s_i\in\Ga(\fU_i,\Oh^*_\cX)$. On the other hand, the $p$-adic norm
$\|s_i(x)\|$ equals $1$ for any element $s_i\in\Ga( \fU_i, \Oh^*_\cX)$
and any $R$-rational point $x\in\fU_i(R)$. Therefore, the $p$-adic
measure on $\fU_i(R)$ that we constructed does not depend on the choice
of a gauge form $\om_i$. Moreover, the $p$-adic measures $\dd\mu_i$ on
$\fU_i(R)$ glue together to a $p$-adic measure $\dd\mu$ on the whole
compact $\cX(R)$, since one has
 \[
 \fU_i(R)\cap\,\fU_j(R)=(\fU_i\cap\,\fU_j)(R)
 \quad\text{for $i,j=1,\dots,k$} 
 \]
and
 \[
 \fU_1(R)\cup\cdots\cup\,\fU_k(R)=
 (\fU_1\cup\cdots\cup\,\fU_k)(R)=\cX(R).
 \]

 \begin{dfn}
 \label{can-m}
 The $p$-adic measure constructed above defined on the set $\cX(R)$ of
$R$-integral points of a $S$-scheme $\cX$ is called the {\em canonical
$p$-adic measure}.
 \end{dfn}

For the canonical $p$-adic measure $\dd\mu$, we obtain the same
property as for the Weil $p$-adic measure $\dd\mu_\om$:

 \begin{thm}
 \label{integ2}
 \[
 \int_{\cX(R)} \dd\mu=\frac{|\cX({F_\fq})|}{q^n}.
 \]
 \end{thm}

 \begin{pf} Using a covering of $\cX$ by some Zariski open subsets
$\fU_1,\dots, \fU_k$, we obtain
 \[
 \int\limits_{\cX(R)}
 \dd\mu=\sum_{i_1}\int\limits_{\fU_{i_1}(R)} \kern-2mm\dd\mu
 -\sum_{i_1<i_2}\kern2mm\int\limits_{(\fU_{i_1}\cap\,\fU_{i_2})(R)}
 \kern-4mm\dd\mu
 \kern2mm+\cdots+\kern2mm
 (-1)^k \kern-4mm\int\limits_{(\fU_1\cap\cdots\cap\,\fU_{k})(R)}
 \kern-4mm\dd\mu
 \]
and
 \begin{multline*}
 \bigl|\cX({F_\fq})\bigr|=
 \sum_{i_1}\bigl|\fU_{i_1}({F_\fq})\bigr|-\sum_{i_1<i_2}\bigl|
(\fU_{i_1}\cap\,\fU_{i_2})({F_\fq})\bigr| \\
 +\cdots+(-1)^k\bigl|(\fU_1\cap\cdots\cap\,\fU_{k})({F_\fq})\bigr|.
 \end{multline*}
It remains to apply Theorem~\ref{weil1} to every intersection
$\fU_{i_1}\cap\cdots\cap\,\fU_{i_s}$.
 \qed\end{pf}

 \begin{thm}
 \label{m-zero}
Let $\cX$ be a regular integral $S$-scheme and $\cZ\subset \cX$ a closed
reduced subscheme of codimension $\ge1$. Then the subset $\cZ(R) \subset
\cX(R)$ has zero measure with respect to the canonical $p$-adic measure
$\dd\mu$ on $\cX(R)$.
 \end{thm}

 \begin{pf} Using a covering of $\cX$ by Zariski open affine subsets
$\fU_1, \dots, \fU_k$, we can always reduce to the case when $\cX$ is
an affine regular integral $S$-scheme and $\cZ\subset\cX$ an
irreducible principal divisor defined by an equation $f=0$, where $f$
is a prime element of $A=\Ga(\cX, \Oh_\cX)$.

Consider the special case $\cX=\A^n_S=\Spec R[X_1,\dots,X_n]$ and
$\cZ=\A^{n-1}_{S} =\Spec R[X_2,\dots,X_n]$, that is, $f=X_1$. For every
positive integer $m$, we denote by $\cZ_m(R)$ the subset in
$\A^n(R)$ consisting of all points $x=(x_1,\dots,x_d)\in R^n $ such
that $x_1\in\fq^m$. One computes the $p$-adic integral in the
straightforward way:
 \[
 \int\limits_{\cZ_m(R)} \bdx=\int\limits_{\A^1(\fq^m)}
 \kern-2mm {\dd x_1}\kern2mm
 \prod_{i=2}^n \left(\int_{\A^1(R)}{\dd x_i}\right)
 =\frac{1}{q^m}.
 \]
On the other hand, we have
 \[
 \cZ(R)=\bigcap_{m =1}^{\infty} \cZ_m(R).
 \]
Hence
 \[
 \int_{\cZ(R)} \bdx=\lim_{m \to\infty }\int_{\cZ_m(R)} \bdx=0,
 \]
and in this case the statement is proved. Using the Noether
normalization theorem reduces the more general case to the above
special one.
 \qed\end{pf}

\section{The Betti numbers}

 \begin{prop}
 \label{main-th2}
Let $X$ and $Y$ be birational smooth projective
$n$-dimensional algebraic varieties over $\C$ having trivial canonical
line bundles. Then there exist Zariski open dense subsets $U\subset
X$ and $V \subset Y$ such that $U$ is isomorphic to $V$ and
$\codim_X(X\setminus U),\codim_Y(Y\setminus V)\ge2$.
 \end{prop}

 \begin{pf}
Consider a birational rational map $\fie\colon X\dasharrow Y$. Since
$X$ is smooth and $Y$ is projective, $\fie$ is regular at the general
point of any prime divisor of $X$, so that there exists a maximal
Zariski open dense subset $U\subset X$ with $\codim_X(X\setminus
U)\ge2$ such that $\fie$ extends to a regular morphism $\fie_0\colon
U\to Y$. Since
$\fie^*\om_Y$ is proportional to $\om_X$, the morphism $\fie_0$ is
\'etale, that is, $\fie_0$ is an open embedding of $U$ into the maximal
open subset $V \subset Y$ where $\fie^{-1}$ is defined. Similarly
$\fie^{-1}$ induces an open embedding of $V$ into $U$, so we conclude
that $\fie_0$ is an isomorphism of $U$ onto $V$.
 \qed\end{pf}
\bigskip

 \begin{pfof}{Theorem~\ref{main-th}} Let $X$ and $Y$ be smooth
projective birational varieties of dimension $n$ over $\C$ with trivial
canonical bundles. By Proposition~\ref{main-th2}, there exist Zariski
open dense subsets $U\subset X$ and $V\subset Y$ with
$\codim_X(X\setminus U)\ge2$ and $\codim_Y(Y\setminus V)\ge2$ and an
isomorphism $\fie\colon U \to V$.

By standard arguments, one can choose a finitely generated
$\Z$-subalgebra $\cR\subset\C$ such that the projective varieties $X$
and $Y$ and the Zariski open subsets $U\subset X$ and $V\subset Y$ are
obtained by base change $*\times_\cS\Spec\C$ from regular projective
schemes $\cX$ and $\cY$ over $\cS:=\Spec\cR$ together with Zariski open
subschemes $\fU\subset\cX$ and $\cV\subset\cY$ over $\cS$. Moreover,
one can choose $\cR$ in such a way that both relative canonical line
bundles $\Om^n_{\cX/\cS}$ and $\Om^n_{\cY/\cS}$ are trivial, both
codimensions $\codim_\cX(\cX\setminus\fU)$ and
$\codim_\cY(\cY\setminus\cV)$ are $\ge2$, and the isomorphism
$\fie\colon U\to V$ is obtained by base change from an isomorphism
$\Phi\colon\fU\to\cV$ over $\cS$.

For almost all prime numbers $p\in\N$, there exist a regular
$R$-integral point $\pi\in\cS \times_{\Spec\Z}{\Spec\Z_p}$, where $R$
is the maximal compact subring in a local $p$-adic field $F$; let $\fq$
be the maximal ideal of $R$. By an appropriate choice of
$\pi\in\cS\times_{\Spec\Z}{\Spec\Z_p}$, we can ensure that both $\cX$
and $\cY$ have good reduction modulo $\fq$. Moreover, we can assume that
the maximal ideal $I(\overline{\pi})$ of the unique closed point in
 \[
 S: =\Spec R \stackrel{\pi}{\hookrightarrow}
\cS\times_{\Spec\Z}{\Spec\Z_p}
 \]
is obtained by base change from some maximal ideal
$J(\overline{\pi})\subset\cR$ lying over the prime ideal $(p)\subset\Z$.

Let $\om_\cX$ and $\om_\cY$ be gauge forms on $\cX$ and $\cY$
respectively and $\om_\fU$ and $\om_\cV$ their restriction to $\fU$
(respectively $\cV$). Since $\Phi^*$ is an isomorphism over $\cS$,
$\Phi^*\om_\cY$ is another gauge form on $\fU$. Hence there exists a
nowhere vanishing regular function $h\in\Ga(\fU,\Oh^*_\cX)$ such that
 \[
 \Phi^* \om_\cV=h \om_\fU.
 \]
The property $\codim_\cX(\cX\setminus \fU)\ge2$ implies that $h$ is
an element of $\Ga(\cX, \Oh^*_\cX)=\cR^*$. Hence, one has $\| h(x) \|
=1$ for all $x\in\cX(F)$, that is, the Weil $p$-adic measures on
$\fU(F)$ associated with $\Phi^* \om_\cV$ and $\om_\fU$ are the same.
The latter implies the following equality of the $p$-adic integrals
 \[
 \int_{\fU(F)} \dd\mu_\cX=\int_{\cV(F)} \dd\mu_\cY.
 \]
By Theorem~\ref{m-zero} and Remark~\ref{point}, (ii), we obtain
 \[
 \int_{\fU(F)} \dd\mu_\cX=\int_{\cX(F)} \dd\mu_\cX =\int_{\cX(R)}
 \dd\mu_\cX
 \]
and
 \[\int_{\cV(\cF)} \dd\mu_\cY=\int_{\cY(\cF)} \dd\mu_\cY=\int_{\cY(R)} \dd\mu_\cY.
 \]
Now, applying the formula in Theorem~\ref{integ2}, we come to the
equality
 \[
 \frac{|\cX({F_\fq})|}{q^n}=\frac{|\cY({F_\fq})|}
{q^n}.
 \]
This shows that the numbers of $F_\fq$-rational points in $\cX$ and
$\cY$ modulo the ideal $J(\overline{\pi})\subset\cR$ are the same. We
now repeat the same argument, replacing $R$ by its cyclotomic extension
$\cR^{(r)}\subset\C$ obtained by adjoining all complex $(q^r-1)$th
roots of unity; we deduce that the projective schemes $\cX$ and $\cY$
have the same number of rational points over $F_\fq^{(r)}$, where
$F_\fq^{(r)}$ is the extension of the finite field $F_\fq$ of degree
$r$. We deduce in particular that the Weil zeta functions
 \[
 Z(\cX,p,t)=\exp \left(\sum_{r =1}^{\infty}
 |\cX({F_\fq^{(r)}})| \frac{t^r}{r} \right)
 \]
and
 \[
 Z(\cY, p, t) =\exp \left(\sum_{r =1}^{\infty}
 |\cY({F_\fq^{(r)}})| \frac{t^r}{r} \right)
 \]
are the same. Using the Weil conjectures proved by Deligne
\cite{Deligne} and the comparison theorem between the \'etale and
singular cohomology, we obtain
 \begin{equation}
 \label{eq_zeta}
 Z(\cX,p,t)=\frac
 {P_1(t)P_3(t)\cdots P_{2n-1}(t)}
 {P_0(t)P_2(t)\cdots P_{2n}(t)}
 \end{equation}
and
 \[
 Z(\cY,p, t)=\frac{ Q_1(t) Q_3(t) \cdots Q_{2n-1}(t)}{
 Q_0(t) Q_2(t) \cdots Q_{2n}(t) },
 \]
where $P_i(t)$ and $Q_i(t)$ are polynomials with integer coefficients
having the properties
 \begin{equation}
 \label{eq_betti}
 \deg P_i(t)=\dim H^i(X, \C), \quad \deg Q_i(t)=\dim H^i(Y, \C)
 \quad \text{for all $i\ge0$.}
 \end{equation}
Since the standard archimedean absolute value of each root of
polynomials $P_i(t)$ and $Q_i(t)$ must be $q^{-i/2}$ and
$P_i(0)=Q_i(0)=1$ for all $i\ge0$, the equality $Z(\cX,p,t)=Z(\cY,p,t)$
implies $P_i(t)=Q_i(t)$ for all $i\ge0$. Therefore, we have $\dim
H^i(X,\C)=\dim H^i(Y,\C)$ for all $i\ge0$.
 \qed\end{pfof}

\section{Further results}

 \begin{dfn} Let $\fie\colon X\dasharrow Y$ be a birational map
between smooth algebraic varieties $X$ and $Y$. We say that $\fie$ {\em
does not change the canonical class}, if for some Hironaka resolution
$\al\colon Z\to X$ of the indeterminacies of $\fie$ the composite
$\al\circ\fie$ extends to a morphism $\be\colon Z\to Y$ such that
$\be^*\Om^n_Y\iso\al^*\Om^n_X$.
 \end{dfn}

The statement of Theorem~\ref{main-th} can be generalized to the case of
birational smooth projective algebraic varieties which do not necessary
have trivial canonical classes as follows:

 \begin{thm}
 \label{main-th3}
Let $X$ and $Y$ be irreducible birational smooth $n$-dimensional
projective algebraic varieties over $\C$. Assume that the exists a
birational rational map $\fie\colon X\dasharrow Y$ which does not
change the canonical class. Then $X$ and $Y$ have the same Betti
numbers.
 \end{thm}

 \begin{pf} We repeat the same arguments as in the proof of
Theorem~\ref{main-th} with the only difference that instead of the Weil
$p$-adic measures associated with gauge forms we consider the canonical
$p$-adic measures (see Definition~\ref{can-m}). Using the birational
morphisms $\al\colon\cZ\to\cX$ and $\be\colon\cZ\to\cY$ having the
property
 \[
 \be^*\Om^n_{\cY/S} \iso \al^*\Om^n_{\cX/S},
 \]
we conclude that for some prime $p\in\N$, the integrals of the canonical
$p$-adic measures $\mu_\cX$ and $\mu_\cY$ over $\cX(R)$ and $\cY(R)$ are
equal, since there exists a dense Zariski open subset $\fU\subset\cZ$ on
which we have $\al^*\mu_\cX=\be^*\mu_\cY$. By Theorem~\ref{integ2}, the
zeta functions of $\cX$ and $\cY$ must be the same.
 \qed\end{pf}

Another immediate application of our method is related to the
McKay correspondence \cite{R}.

 \begin{thm}
Let $G \subset\SL(n,\C)$ be a finite subgroup. Assume that there
exist two different resolutions of singularities on $W:=\C^n/G$:
 \[
 f\colon X \to W, \quad g\colon Y \to W
 \]
such that both canonical line bundles $\Om^n_X$ and $\Om^n_Y$
are trivial. Then the Euler numbers of $X$ and $Y$ are the same.
 \end{thm}

 \begin{pf} We extend the varieties $X$ and $Y$ to regular schemes over
a scheme $\cS$ of finite type over $\Spec\Z$. Moreover, one can choose
$\cS$ in such a way that the birational morphisms $f$ and $g$ extend to
birational $\cS$-morphisms
 \[
 F\colon \cX \to\cW, \quad G\colon \cY \to\cW,
 \]
where $\cW$ is a scheme over $\cS$ extending $W$. Using the same
arguments as in the proof of Theorem~\ref{main-th}, one obtains that
there exists a prime $p\in\N$ such that $Z(\cX,p,t)=Z(\cY,p,t)$. On the
other hand, in view of (\ref{eq_betti}), the Euler number is determined
by the Weil zeta function (\ref{eq_zeta}) as the degree of the numerator
minus the degree of the denominator. Hence $e(X)=e(Y)$.
 \qed\end{pf}

With a little bit more work one can prove even more precise
statement:

 \begin{thm} Let $G\subset\SL(n,\C)$ be a finite subgroup and
$W:=\C^n/G$. Assume that there exists a resolution
 \[
 f\colon X \to W
 \]
with trivial canonical line bundle $\Om^n_X$. Then the Euler
number of $X$ equals the number of conjugacy classes in $G$.
 \end{thm}

 \begin{rem} As we saw in the proof of Theorem~\ref{main-th2}, the Weil
zeta functions of $Z(\cX,p,t)$ and $Z(\cY,p,t)$ are equal for almost
all primes $p\in\Spec\Z$. This fact being expressed in terms of the
associated $L$-functions indicates that the isomorphism $H^i(X,\C)\iso
H^i(Y,\C)$ for all $i\ge0$ we have established must have some more deep
motivic nature. Recently Kontsevich suggested an idea of a motivic
integration \cite{K}, developed by Denef and Loeser \cite{DL}. In
particular, this technique allows to prove that not only the Betti
numbers, but also the Hodge numbers of $X$ and $Y$ in \ref{main-th}
must be the same.
 \end{rem}

\noindent
Victor V. Batyrev \\
Mathematisches Institut, Universit\"at T\"ubingen \\
Auf der Morgenstelle 10, 72076 T\"ubingen, Germany \\
e-mail: batyrev@bastau.mathematik.uni-tuebingen.de

 \end{document}